\documentclass[10pt]{amsart}
\usepackage{amsmath}

\newtheorem{thm}{Theorem}[section]
\newtheorem{lem}[thm]{Lemma}

\newtheorem{cor}[thm]{Corollary}
\newtheorem{defn}[thm]{Definition}
\newtheorem*{remark}{Remark}

\numberwithin{equation}{section}

\newcommand{\R}{{\mathbb R}}
\newcommand{\C}{{\mathbb C}}
\newcommand{\N}{{\mathbb N}}

\newcommand{\supp}{\operatorname{supp}}
\renewcommand{\mod}{\;\operatorname{mod}}
\newcommand{\Tr}{\operatorname{Tr}}

\newcommand{\intinf}{\int_{-\infty}^\infty}
\renewcommand{\idotsint}{\int\!\cdots\int}
\newcommand{\nidotsint}[1]{\underbrace{\int\!\cdots\;}_{#1}\!\!\!\int}

\newcommand{\ave}[1]{\left\langle#1\right\rangle} 
\newcommand{\aveq}[1]{\left\langle#1\right\rangle_q} 
\newcommand{\E}{{\mathbb{E}}} 

\newcommand{\Wf}{W_f} 
\newcommand{\Zf}{Z_f} 

\newcommand{\I}{1\!\!1} 
\renewcommand{\Re}{{\mathfrak{Re}}}
\renewcommand{\Im}{{\mathfrak{Im}}}
\renewcommand{\O}{{\mathcal{O}}} 

\renewcommand{\a}{\alpha}
\renewcommand{\b}{\beta}
\newcommand{\g}{\gamma}
\renewcommand{\l}{\lambda}
\renewcommand{\t}{\theta}

\renewcommand{\i}{{\mathrm{i}}} 
\renewcommand{\d}{{\mathrm{d}}} 
\renewcommand{\^}{\widehat} 

\begin{document}
\title{Linear statistics of low-lying zeros of $L$--functions}
\author{C.P. Hughes and Z. Rudnick}
\address{Raymond and Beverly Sackler School of Mathematical Sciences,
Tel Aviv University, Tel Aviv 69978, Israel. Current address:
American Institute of Mathematics, 360 Portage Ave, Palo Alto, CA
94306-2244, USA ({\tt hughes@aimath.org})}
\address{Raymond and Beverly Sackler School of Mathematical Sciences,
Tel Aviv University, Tel Aviv 69978, Israel ({\tt
rudnick@post.tau.ac.il})}

\date{3 December 2002}

\thanks{Supported in part by  the EC TMR network
``Mathematical aspects of Quantum Chaos'', EC-contract no
HPRN-CT-2000-00103.}

\begin{abstract}
We consider linear statistics of the scaled zeros of Dirichlet
$L$--functions, and show that the first few moments converge to
the Gaussian moments. The number of Gaussian moments depends on
the particular statistic considered. The same phenomenon is found
in Random Matrix Theory, where we consider linear statistics of
scaled eigenphases for matrices in the unitary group. In that case
the higher moments are no longer Gaussian. We conjecture that this
also happens for Dirichlet $L$--functions.
\end{abstract}

\maketitle

\section{Introduction}

Let $q$ be an odd prime and $\chi$ a Dirichlet character modulo
$q$. For $\Re(s)>1$ the Dirichlet $L$--function $L(s,\chi)$ is
defined as
$$
L(s,\chi):=\sum_{n=1}^\infty \frac {\chi(n)}{n^s} = \prod_p
\left(1-\frac{\chi(p)}{p^s}\right)^{-1}
$$

Each function $L(s,\chi)$ has an infinite set of non-trivial zeros
$1/2+\i\g_{\chi,j}$ which can be ordered so that
$$
\dots \leq \Re(\g_{\chi,-2}) \leq \Re(\g_{\chi,-1}) < 0 \leq
\Re(\g_{\chi,1}) \leq \Re(\g_{\chi,2}) \leq \dots
$$
Note that we don't assume the Generalised Riemann Hypothesis (GRH)
since we allow the $\g_{\chi,j}$ to be complex.

Denote by $N(T,\chi)$ number of non-trivial zeros such that
$0<\Re(\g_{\chi,j})< T$. Then  for fixed $T>0$,
$$
\frac{1}{q-2} \sum_{\chi\neq\chi_0} N(T,\chi)\sim \frac{T}{2\pi}
\log\frac{q T}{2\pi e}  \quad \text{ as }q\to\infty
$$
where the sum is taken over all the $q-2$ non-trivial characters
modulo the prime $q$, see Titchmarsh \cite{Titchmarsh}, Siegel
\cite{Siegel},  Selberg \cite{selberg_L}. We will therefore scale
the zeros by defining
\begin{equation*}
x_{\chi,j} := \frac{\log q}{2\pi} \g_{\chi,j} .
\end{equation*}

The purpose of this paper is to consider linear statistics of the
low-lying $x_{\chi,j}$. Let $f$ be a rapidly decaying even test
function, and consider the linear statistic
\begin{equation*}
\Wf(\chi) := \sum_{j=-\infty}^\infty f(x_{\chi,j})
\end{equation*}

Linear statistics for low-lying zeros of several families of
$L$-functions were investigated systematically by Katz and Sarnak
\cite{KatzSarnak} and  by Iwaniec, Luo and Sarnak \cite{ILS} where
they were called ``one-level densities''. We prefer to use the
terminology ``linear statistic'' which is traditional in random
matrix theory.

Define the $\chi$--average of $\Wf(\chi)$ as
\begin{equation*}
\ave{\Wf}_q := \frac{1}{q-2} \sum_{\chi\neq\chi_0} \Wf(\chi)
\end{equation*}
In order to understand the distribution of $\Wf(\chi)$ we
calculate its first few moments $\ave{\Wf^m}_q$ . In
\S\ref{sect:W_mean} we prove that if $\supp \^f \subseteq [-2,2]$,
then the mean of $\Wf(\chi)$ is $\int_{-\infty}^\infty f(x)\;\d x
$, and in \S\ref{sect:variance of Wf} we show that the variance
converges to $\int_{-1}^1 |u| \^f(u)^2\;\d u$ if $\supp \^f
\subseteq[-1,1]$. In \S\ref{sect:mmts_W} we show that if $\supp
\^f \subset (-2/m,2/m)$ (for $m\geq 2$) then the first $m$ moments
of $\Wf$ converge to the first $m$ moments of a normal random
variable with mean and variance as above.

If all moments of $\Wf(\chi)$ were Gaussian, then we would be able
to conclude that $\Wf$ was normally distributed. Indeed, it
follows from the work of Selberg \cite{selberg_L} that scaling the
$\g_{\chi,j}$ by anything much less than $\log q$, leads to a
Gaussian distribution. However, with scaling on the order of $\log
q$, this cannot be the case for all $f$ since taking $f$ to be an
indicator function, the limiting distribution is discrete. We
therefore say that $\Wf$ displays {\em mock-Gaussian} behaviour.

As suggested by Katz and Sarnak \cite{KatzSarnak}, one may try to
model properties of $\Wf(\chi)$ by  random matrix theory. Let $U$
be an $N\times N$ unitary matrix, with eigenvalues $e^{\i\t_n}$.
The statistical distribution of $\frac{N}{2\pi} \t_n$ has been
conjectured to converge to the empirical distribution of
$x_{\chi,j}$ as $q$ and $N$ both tend to infinity

Therefore, as an aide to understanding
$\Wf(\chi)=\sum_{j=-\infty}^\infty f(x_{\chi,j})$, one might wish
to calculate the moments of $\sum_{n=1}^N f(\frac{N}{2\pi}\t_n)$.
However, since the $\t_n$ are angles, it is more natural (and
indeed more convenient) to consider the $2\pi$--periodic function
\begin{equation*}
F_N(\t) := \sum_{j=-\infty}^\infty f\left(\tfrac{N}{2\pi} (\t+2\pi
j)\right)
\end{equation*}
and model $\Wf(\chi)$ by
\begin{equation*}
\Zf(U) := \sum_{j=1}^N F_N(\t_j)
\end{equation*}
where  $U$ is an $N\times N$ unitary matrix with eigenangles
$\t_1,\dots,\t_N$. Note that the scaling $N/2\pi$ (the mean
density) is equivalent to the scaling $\frac{\log q}{2\pi}$ for
the zeros of $L$--functions.

Our results for $\Zf(U)$ are given in
\S\S\ref{sect:rmt_model},\ref{sect:Z_f_not_fourier}. Writing $\E$
to denote the average over the unitary group with Haar measure,
then without any restrictions on the support of the function $f$,
we prove that $\E\{\Zf(U)\} = \int_{-\infty}^\infty f(x)\;\d x$,
and that the variance tends to
\begin{equation}\label{eq:rmt_var}
\sigma(f)^2 = \int_{-\infty}^\infty \min(1,|u|) \^f(u)^2\;\d u
\end{equation}
Observe that this is in complete agreement with the mean and
variance of $\Wf(\chi)$ if $\^f$ has the same support
restrictions. Furthermore, we show in
\S\ref{sect:Zf_mock_Gaussian} that for any integer $m\geq 2$, if
$\supp \^f \subseteq [-2/m,2/m]$, then
\begin{equation*}
\lim_{N\to\infty} \E \left\{ \left( \Zf-\E\{\Zf\} \right)^m
\right\} =
\begin{cases}
0 & \text{ if } m \text{ odd}\\
\frac{m!}{2^{m/2}  (m/2)!} \sigma^m & \text{ if $m$ even}
\end{cases}
\end{equation*}
where $\sigma^2$, the variance, is given in \eqref{eq:rmt_var}.
These are the moments of a normal random variable, so again we see
mock-Gaussian behaviour, with the same restrictions on the support
of $\^f$ as in $\Wf(\chi)$.

To understand the mock-Gaussian behaviour, note that if we had
defined
$$
F_N^{(L)}(\t) =  \sum_{j=-\infty}^\infty f\left(L (\t+2\pi
j)\right)
$$
where $L\to\infty$ subject to $\frac{L}{N}\to 0$, then Soshnikov
\cite{Sosh} (see also \cite{DE}) has shown that the mean of
$\Zf^{(L)}(U)$ converges to $\frac{N}{2\pi L}\int_{-\infty}^\infty
f(x)\;\d x$, and that the centered random variable
$\Zf^{(L)}-\E\left\{\Zf^{(L)}\right\}$ converges in distribution
to a normal random variable with mean zero and variance
$\int_{-\infty}^\infty \^f(u)^2 |u|\;\d u$. Our scaling is
$L=\frac{N}{2\pi}$, which is just outside the range of Soshnikov's
result. Indeed, note that the variance \eqref{eq:rmt_var} is
different if $\supp \^f \not\subseteq [-1,1]$.

In \S\ref{sect:Z_f_not_fourier} we show that all moments of $\Zf$
can be calculated exactly within random matrix theory, without any
restrictions on the support of $\^f$. They are given by a
complicated expression, but are certainly not Gaussian moments in
general. The moments of $\Zf(U)$ grow sufficiently slowly that
they uniquely determine its distribution.

Finally, in \S\ref{sect:small_zeros} we apply the results of
\S\S\ref{sect:W_mean},\ref{sect:variance of Wf} to show that,
under the assumption of GRH, for each $q$, there exist $\chi$ such
that the height of the lowest zero is less than $1/4$ times the
expected height. We also obtain a similar result where a positive
proportion of $L$--functions have their first zero less than
$0.633$ times the expected height.

Linear statistics of the high zeros of a \emph{fixed}
$L$--function also show mock-Gaussian behaviour \cite{hr3}, having
the same moments as the linear statistics of low-lying zeros
considered in this paper.

Moments of linear statistics in other classical compact groups,
like $SO(2N)$, $SO(2N+1)$ and $Sp(2N)$ also show mock-Gaussian
behaviour, \cite{hr2}. Other $L$--functions can be modeled by
these groups. For example, in the case of quadratic $L$--functions
mock-Gaussian behaviour can be deduced from the work of Rubinstein
\cite{Rubin}. Specifically, if $\supp \^f \subset
(-\frac{1}{m},\frac{1}{m})$ then the first $m$ moments are
Gaussian with mean $\^f(0)-\int_0^1 \^f(u)\;\d u$ and variance
$4\int_0^{1/2} u \^f(u)^2\;\d u$, exactly as in the group
$Sp(2N)$. We remark that this is half the unitary range. Assuming
GRH, the results of \"Ozl\"uck and Snyder \cite{OzSn} show that
the mean is indeed $\^f(0)-\int_0^1 \^f(u)\;\d u$ so long as
$\supp \^f \subset (-2,2)$, in the sense that
$$\frac{1}{D}\sum_{d} e^{-\pi d^2/D^2} \sum_{\g} f(\frac{\g \log
D}{2\pi}) = \^f(0) - \int_0^1 \^f(u)\;\d u + o(1)$$

We note that the arguments given here to study moments of the
linear statistic show that the ``$n$-level densities''
\cite{KatzSarnak, Rubin} of this family of $L$-functions coincide
with those of the unitary group, in a suitable range, by purely
combinatorial arguments. Since that is not our purpose here, we
leave it for the reader.

Throughout all this paper, the Fourier transform is $\^ f(u) =
\int_{-\infty}^\infty f(x) e^{-2\pi\i xu}\;\d x$, and thus the
inverse transform is $f(x) = \int_{-\infty}^\infty \^f(u)
e^{2\pi\i xu}\;\d u$.

\section{The scaled level density $\Wf$}

\subsection{Zeros and the Explicit Formula}

A Dirichlet character $\chi:\N\longrightarrow \C$ is a function
such that $\chi(n+q)=\chi(n)$ for all $n$; $\chi(n)=0$ if $n$ and
$q$ have a common divisor; and $\chi(mn)=\chi(m)\chi(n)$ for all
$m,n$. We say $\chi$ is the trivial character modulo $q$ (denoted
$\chi_0$) if $\chi(n)=1$ for all $n$ coprime to $q$. We say that
$\chi$ is even if $\chi(-1)=1$, and odd if $\chi(-1)=-1$.

The Explicit Formula is the following relation between a sum over
zeros of $L(s,\chi)$ and a sum over prime powers. To describe it,
let
$$
a(\chi) = \begin{cases} 0,& \chi \text{ even}\\
                       1,& \chi \text{ odd}
\end{cases}
$$
and let $h(r)$ be any even analytic function in the strip $-c\leq
\Im(r)\leq 1+c$ (for $c>0$) such that $|h(r)|\leq
A(1+|r|)^{-(1+\delta)}$ (for $r\in\R$, $A>0$, $\delta>0$). Set
$g(u) = \frac 1{2\pi}\int_{-\infty}^\infty h(r)e^{-\i ru} \;\d r$,
so that $h(r) = \intinf g(u) e^{\i ru} \;\d u$. Then
\begin{multline}\label{Explicit formula for L(s,chi)}
\sum_j h(\gamma_{j,\chi}) =
\frac 1{2\pi}\intinf h(r) \left(\log q+ G_\chi (r)\right) \;\d r \\
-  \sum_n \frac{\Lambda(n)}{\sqrt{n} }g(\log n)
\left(\chi(n)+\bar\chi(n)\right)
\end{multline}
where
$$G_\chi(r)=\frac{\Gamma'}{\Gamma}(\frac 12+a(\chi)+\i r)+
\frac{\Gamma'}{\Gamma}(\frac 12+a(\chi)-\i r)
-\frac{1}{2}\log\pi\;.
$$
and the von Mangoldt function $\Lambda(n)$ is defined as $\log p$
if $n=p^k$ is a prime power, and zero otherwise.

\subsection{A decomposition of $\Wf$}
For test functions $f$ define the scaled level density $\Wf(\chi)$
as
$$
\Wf(\chi):=\sum_{j} f(\frac{\log q}{2\pi} \gamma_{\chi,j})
$$
the sum over all nontrivial zeros of $L(s,\chi)$.

\begin{defn}
$f(x)$ is an admissible test functions for $\Wf(\chi)$ if it is a
real, even function, whose Fourier transform $\widehat f(u)$ is
compactly supported, and such that $f(r) \ll (1+|r|)^{-1-\delta}$.
\end{defn}

We will transform $\Wf$ into a sum over prime powers by using the
explicit formula  for $L(s,\chi)$. In \eqref{Explicit formula for
L(s,chi)}, take $h(r)=f(\frac {\log q}{2\pi} r)$, so that $g(u) =
\frac 1{\log q}\widehat f(\frac u{\log q})$, and note that the
conditions on $\^f$ easily imply the analyticity and decay
condition on $h(r)$ in the explicit formula. We then get a
decomposition of $\Wf(\chi)$ as
\begin{equation}\label{Decomposition of \Wf}
\Wf(\chi) = \overline{\Wf}(\chi)  + \Wf^{osc}(\chi)
\end{equation}
where
\begin{equation*}
\overline{\Wf}(\chi)   := \frac 1{2\pi}  \intinf f(\frac {\log
q}{2\pi} r) \left(\log q + G_\chi(r) \right) \;\d r
\end{equation*}
and  an oscillatory term
\begin{equation}\label{expression for \Wf^{osc}}
\Wf^{osc}(\chi) := -  \frac 1{\log q} \sum_n
\frac{\Lambda(n)}{\sqrt{n} }\widehat f(\frac {\log n}{\log q})
\left(\chi(n)+\bar\chi(n)\right) \;.
\end{equation}
The first term $\overline{\Wf}(\chi)$   gives
$$
\overline{\Wf}(\chi)= \intinf f(x)\;\d x +\O(\frac 1{\log q})
$$
which is asymptotically independent of $\chi$.

\section{The expectation of $\Wf$} \label{sect:W_mean}
The ``expectation'' of $\Wf$ is defined as the average over all
$q-2$ nontrivial characters modulo prime $q$
$$
\aveq{\Wf}:=\frac 1{q-2} \sum_{\chi\neq \chi_0 } \Wf(\chi)
$$

\begin{thm}\label{thm:expectation}
Let $f$ be an admissible function, and assume $\supp(\widehat f)
\subseteq [-2,2]$. Then as $q\to \infty$,
$$
\aveq{\Wf} = \int_{-\infty}^\infty f(x)\;\d x
+\O\left(\frac{1}{\log q}\right)\; .
$$
\end{thm}
\begin{proof}
We will use the decomposition \eqref{Decomposition of \Wf} and
average over $\chi$:
$$
\aveq{ \Wf}  = \aveq{\overline{\Wf}} +\aveq{\Wf^{osc} } \;.
$$
Since $\overline{\Wf}$ is asymptotically constant for any $f$, we
have
$$
\aveq{\overline{\Wf}}= \intinf f(x)\;\d x +\O(\frac 1{\log q})
$$
and thus it will suffice to show that
$$
\aveq{\Wf^{osc} } =\O(\frac{1}{\log q}) \;.
$$
We will show this under the assumption $\supp(\widehat f)
\subseteq [-2,2]$.

We have by \eqref{expression for \Wf^{osc}} that
\begin{equation*}
\aveq{\Wf^{osc} }  = - \frac 1{\log q} \sum_{n}
\frac{\Lambda(n)}{\sqrt{n} }\widehat f(\frac {\log n}{\log q})
\left(\aveq{\chi(n)}+\aveq{\bar\chi(n)}\right) \;.
\end{equation*}
The mean value of $\chi$ is
$$
\aveq{\chi(n)}=\aveq{ \bar\chi(n) } =
\begin{cases}
1,& n\equiv 1 \mod q\\
0,& q\mid n\\
-\frac 1{q-2},& n\not\equiv 0,1 \mod q
\end{cases}
$$

Thus we find
\begin{multline*}
\aveq{\Wf^{osc}} = \frac {-2}{\log q}\sum_{n\equiv 1\mod q}
\frac{\Lambda(n)}{\sqrt{n}}\widehat f(\frac {\log n}{\log q}) +
\frac{2}{\log q} \frac 1{q-2} \sum_{n\not\equiv 1,0 \mod q}
\frac{\Lambda(n)}{\sqrt{n}}\widehat f(\frac {\log n}{\log q})
\end{multline*}
Assume that $\supp(\widehat f) \subseteq [-\alpha,\alpha]$ for
$\alpha>0$. Then the sum is over $n\leq q^\alpha$. Since $\widehat
f$ is bounded, we may replace it by $1$ over that range, and we
therefore have
\begin{equation}\label{eq:Wf_osc}
\aveq{\Wf^{osc}} \ll \frac 1{\log q}\sum_{\substack{n\equiv 1\mod
q\\n\leq q^\alpha}}  \frac{\Lambda(n)}{\sqrt{n} } + \frac{1}{\log
q} \frac 1{q-2} \sum_{\substack{n\leq q^\alpha}}
\frac{\Lambda(n)}{\sqrt{n} }
\end{equation}

To deal with the first sum in \eqref{eq:Wf_osc} one could replace
primes by integers by noting that $\Lambda(n) \ll \log n$,
obtaining
\begin{align*}
\frac 1{\log q}\sum_{\substack{n\equiv 1\mod q\\n\leq q^\alpha}}
\frac{\Lambda(n)}{\sqrt{n} }
&\ll \frac{1}{\log q} \sum_{\substack{n\equiv 1\mod q\\n\leq q^{\alpha}}} \frac{\log n}{\sqrt{n}}\\
&= \frac{1}{\log q} \sum_{m< q^{\alpha-1}} \frac{\log (mq+1)}{\sqrt{mq+1}}\\
&\ll \frac{1}{\sqrt{q}} q^{(\a-1)/2}
\end{align*}
which vanishes for $\a<2$. However, one gets a slightly stronger
result by the Brun-Titchmarsh Theorem, \cite{MontVau}, which says
that if $\pi(x;q,a)$ is the number of primes $p\leq x$, $p \equiv
a \mod q$, where $q$ and $a$ are coprime, then for $x>2q$,
$$
\pi(x;q,a) < \frac{2 x}{\varphi(q)\log(x/q)}
$$
Therefore, for $\alpha>1$,
\begin{align*}
\frac 1{\log q}\sum_{\substack{n\equiv 1\mod q\\n\leq q^\alpha}}
\frac{\Lambda(n)}{\sqrt{n} }
&\ll \frac 1{\log q}\sum_{\substack{p\equiv 1\mod q\\p\leq q^\alpha}}  \frac{\log p}{\sqrt{p} }\\
&\ll \frac{1}{\log q} \int_{2q}^{q^\a} \frac{\log x}{\sqrt{x}} \frac{1}{q} \frac{\d x}{\log (x/q)}\\
&\ll \frac{1}{\log q} q^{-1+\a/2}
\end{align*}
which vanishes for $\alpha\leq 2$.

To deal with the second sum in \eqref{eq:Wf_osc}, one could
similarly replace primes by integers, and note that
$$
\frac{1}{\log q} \frac 1{q-2} \sum_{\substack{n\leq q^\alpha}}
\frac{\log n}{\sqrt{n} } \ll q^{-1+\a/2}
$$
which vanishes if $\alpha<2$. Again this result can be
strengthened by using the Prime Number Theorem, since
\begin{align*}
\frac{1}{\log q} \frac 1{q-2} \sum_{\substack{n\leq
q^\alpha\\n\not\equiv 0,1 \mod q}} \frac{\Lambda(n)}{\sqrt{n} }
&\ll \frac{1}{\log q} \frac 1{q} \sum_{p\leq q^\alpha} \frac{\log p}{\sqrt{p} }\\
&\ll \frac{1}{\log q} \frac 1{q} \int_2^{q^\alpha} \frac{\log x}{\sqrt{x} } \frac{\d x}{\log x}\\
&\ll \frac{1}{\log q} q^{-1+\alpha/2}
\end{align*}
which vanishes for $\a\leq 2$.

Thus $\aveq{\Wf^{osc}}\ll \frac{q^{-1+\alpha/2}}{\log q}$ and so
if $\alpha\leq 2$ we find $\aveq{\Wf^{osc}} \to 0$ as required.
\end{proof}

\begin{remark}
Set
$$\Wf^{(t)}(\chi) = \sum_{\g_{\chi,j}} f\left(\frac{\log q}{2\pi}(\g_{\chi,j}-t)\right)$$
which is like $\Wf(\chi)$ but with the zeros shifted by height
$t$. Averaging over all characters modulo $q$, and doing an extra
smooth average over $t$, the expected value of $\Wf^{(t)}$
converges to $\intinf f(x)\;\d x$ without any restriction on the
support of $\^f$.
\end{remark}

\section{The variance of $\Wf$}\label{sect:variance of Wf}

\begin{thm}\label{thm:Wf_variance}
Let $f$ be an admissible function and assume $\supp \^f \subseteq
[-1,1]$, then the variance of $\Wf$ tends to
\begin{equation}\label{eq:Wf_variance}
\sigma(f)^2 = \int_{-1}^1 |u| \^ f(u)^2 \;\d u
\end{equation}
\end{thm}

\begin{proof}
The variance of $\Wf$ is, by Theorem \ref{thm:expectation}
\begin{equation*}
\aveq{\left(\Wf - \aveq{\Wf}\right)^2} = \aveq{(\Wf^{osc} )^2}
+\O\left(\frac{1}{\log q}\right)
\end{equation*}
and by \eqref{expression for \Wf^{osc}}
\begin{multline*}
\aveq{(\Wf^{osc} )^2} = \frac{1}{(\log q)^2} \sum_{n_1}\sum_{n_2} \frac{\Lambda(n_1)}{\sqrt{n_1}}  \frac{\Lambda(n_2)}{\sqrt{n_2}} \^f\left(\frac{\log n_1}{\log q}\right)  \^f\left(\frac{\log n_2}{\log q}\right) \times\\
\times
\left(\aveq{\chi(n_1)\chi(n_2)}+\aveq{\chi(n_1)\bar\chi(n_2)}+\aveq{\bar\chi(n_1)\chi(n_2)}+\aveq{\bar\chi(n_1)\bar\chi(n_2)}\right)
\end{multline*}
Now,
$$
\aveq{\chi(n_1)\chi(n_2)} = \aveq{\chi(n_1 n_2)} =
\begin{cases}
1, & n_1 n_2 \equiv 1 \mod q\\
0, & n_1 \text{ or } n_2 \equiv 0 \mod q\\
\frac{-1}{q-2}, & \text{ otherwise}
\end{cases}
$$
and
$$
\aveq{\chi(n_1)\bar\chi(n_2)} =
\begin{cases}
1, & n_1 \equiv n_2 \not\equiv 0 \mod q\\
0, & n_1 \text{ or } n_2 \equiv 0 \mod q\\
\frac{-1}{q-2}, & \text{ otherwise}
\end{cases}
$$
Since we assume $\supp \^f \subseteq [-1,1]$, we need only
consider $n_1, n_2 \leq q$. Therefore, writing $\bar n_1$ for the
inverse of $n_1$ modulo $q$,
\begin{multline}\label{eq:Wf_var_mainterm}
\aveq{(\Wf^{osc})^2} = \frac{1}{(\log q)^2} \Biggl(2\sum_{n_1=2}^{q-1} \frac{\Lambda(n_1)^2}{n_1} \^f\left(\frac{\log n_1}{\log q}\right)^2 \\
+ 2\sum_{n_1=2}^{q-1} \frac{\Lambda(n_1)}{\sqrt{n_1}}  \frac{\Lambda(\bar n_1)}{\sqrt{\bar n_1}} \^f\left(\frac{\log n_1}{\log q}\right)  \^f\left(\frac{\log \bar n_1}{\log q}\right) \Biggr)\\
+ \O\left(\frac{1}{q-2} \left(\frac{1}{\log q}\sum_{n\leq q}
\frac{\Lambda(n)}{\sqrt{n}} \right)^2\right)
\end{multline}
By the Prime Number Theorem,
$$
\frac{1}{\log q}\sum_{n\leq q} \frac{\Lambda(n)}{\sqrt{n}} \ll
\frac{\sqrt{q}}{\log q}
$$
and so the $\O$ term in \eqref{eq:Wf_var_mainterm} is bounded by
$\frac{1}{(\log q)^2}$.

Now,
\begin{align*}
\frac{2}{(\log q)^2}
\sum_{n_1=2}^{q-1}\frac{\Lambda(n_1)}{\sqrt{n_1}}
\frac{\Lambda(\bar n_1)}{\sqrt{\bar n_1}} \^f\left(\frac{\log
n_1}{\log q}\right)  \^f\left(\frac{\log \bar n_1}{\log q}\right)
&\ll \sum_{\substack{p<q\\p,\bar p \text{  prime}}} \frac{1}{\sqrt{p}}\frac{1}{\sqrt{\bar p}}\\
&= \sum_{\substack{p<q\\p,\bar p \text{  prime}}}
\frac{1}{\sqrt{k_p q+1}}
\end{align*}
where $k_p$ is defined so that $p \bar p = 1+k_p q$. By unique
factorization, there are exactly two primes less than $q$ which
produce a given $k_p$, namely $p$ and $\bar p$. Therefore there
can be at most $\tfrac{1}{2}\pi(q)$ different $k_p$, all lying
between $1$ and $q-1$, and so
\begin{align*}
\frac{1}{(\log q)^2}
\sum_{n_1=2}^{q-1}\frac{\Lambda(n_1)}{\sqrt{n_1}}
\frac{\Lambda(\bar n_1)}{\sqrt{\bar n_1}}
&\ll \sum_{k=1}^{\pi(q)/2} \frac{1}{\sqrt{kq+1}}\\
&\ll \frac{1}{\sqrt{q}}\sqrt{\frac{q}{\log q}}
\end{align*}
Inserting this into \eqref{eq:Wf_var_mainterm} we have
\begin{align*}
\aveq{(\Wf^{osc})^2} &= \frac{2}{(\log q)^2} \sum_{n=2}^{q-1} \frac{\Lambda(n)^2}{n_1} \^f(\frac{\log n}{\log q})^2 + \O\left(\frac{1}{\sqrt{\log q}}\right)\\
&=\frac{2}{(\log q)^2} \sum_{p<q} \frac{(\log p)^2}{p} \^f(\frac{\log p}{\log q})^2 + \O\left(\frac{1}{\sqrt{\log q}}\right)\\
&= \frac{2}{(\log q)^2} \int_2^q \frac{(\log x)^2}{x} \^f(\frac{\log x}{\log q})^2 \frac{\d x}{\log x} + \O\left(\frac{1}{\sqrt{\log q}}\right)\\
&= 2\int_{\log 2/\log q}^1 u \^f(u)^2\;\d u +
\O\left(\frac{1}{\sqrt{\log q}}\right)
\end{align*}
and so we see that, since $\^f(u)$ is an even function,
$$
\lim_{q\to\infty}\aveq{(\Wf^{osc})^2} = \int_{-1}^1 |u|
\^f(u)^2\;\d u
$$
as required.
\end{proof}

\begin{remark}
Assuming GRH, \"Ozl\"uck \cite{Oz} shows that the variance of
linear statistics of the scaled zeros shifted by $t$ and weighted
by a smooth function $K$, converge to the weighted form of
\eqref{eq:Wf_variance} so long as $\supp \^f\subset (-2,2)$, when
averaged over $t$ and over all characters of modulus less than
$q$. In fact, random matrix theory suggests that
\eqref{eq:Wf_variance} is the correct variance for all admissible
functions (Theorem~\ref{thm:Zf_var}).
\end{remark}

\section{The moments of $\Wf$} \label{sect:mmts_W}

We now attempt to understand the {\em distribution} of the scaled
level density $\Wf$ around its expected value. We will find that
the first few moments of $\Wf$ converge to those of a Gaussian
random variable with mean $\lim \overline{\Wf}
=\int_{-\infty}^\infty f(x)\;\d x$ and variance
$$
\int_{-\infty}^\infty \min(1,|u|) \^ f(u)^2 \;\d u \; .
$$

\begin{thm} \label{thm-L:Gaussian moments}
Let $f$ be an admissible function, and assume that
\begin{equation} \label{support condition}
\supp \widehat f \subseteq [-\alpha,\alpha],\quad \alpha>0\;.
\end{equation}
If $m<2/\alpha$, then the $m$-th  moment of $\Wf^{osc}$ is
$$
\lim_{q\to \infty} \aveq{ ( \Wf^{osc} )^m }  =
\begin{cases}
\frac{m!}{2^{m/2} (m/2)!} \sigma(f)^{m},& m \text{ even}\\
0,& m\text{ odd}
\end{cases}
$$
where $\sigma(f)^2$, the variance, is given in
\eqref{eq:Wf_variance}.
\end{thm}

\begin{proof}
By \eqref{expression for \Wf^{osc}}, we have
\begin{equation*}
\Wf^{osc}(\chi)  = -\frac 1{\log q} \sum_n
\frac{\Lambda(n)}{\sqrt{n}} \widehat f(\frac{\log n}{\log q})
\left ( \chi(n)+\bar\chi(n) \right)
\end{equation*}

This gives
\begin{equation*}
\begin{split}
(\Wf^{osc}(\chi) )^m & = (-\frac 1{\log q})^m \prod_{j=1}^m
\sum_{n_j} \frac{\Lambda(n_j)}{\sqrt{n_j}} \widehat f(\frac{\log
n_j}{\log q})
\left( \chi(n_j)+\bar\chi(n_j) \right) \\
&=\sum_{S\subset\{1,\dots,m\}} J(S)(\chi)
\end{split}
\end{equation*}
where for any subset of indices $S\subseteq \{1,\dots,m\}$, the
summand $J(S)(\chi)$ corresponds to the different ways of picking
$\chi$ and $\bar\chi$:
$$
J(S)(\chi)  := (-\frac 1{\log q})^m \sum_{n_1,\dots ,n_m}
\prod_{j=1}^m \frac{\Lambda(n_j)}{\sqrt{n_j}} \widehat
f(\frac{\log n_j}{\log q}) \chi(\prod_{j\in S}n_j)
\bar\chi(\prod_{i\notin S}n_i) \;.
$$

Now average over the nontrivial characters $\chi$, using
$$
\aveq{  \chi(\prod_{j\in S}n_j) \bar\chi(\prod_{i\notin S}n_i) } =
\begin{cases}
1,& \prod_{j\in S}n_j\equiv \prod_{i\notin S}n_i \not\equiv 0 \mod q \\
0,& \text{ any one of } n_j \equiv 0 \mod q\\
-1/(q-2),& \text{ otherwise}
\end{cases}
$$
This gives
$$
\aveq{(\Wf^{osc}(\chi) )^m} = \sum_{S\subset\{1,\dots,m\}}
\aveq{J(S)}
$$
with
\begin{multline}\label{J(S)+O}
\aveq{ J(S) } = (-\frac 1{\log q})^m \sum_{\prod_{j\in S}n_j\equiv
\prod_{i\notin S}n_i \mod q} \prod_{j=1}^m
\frac{\Lambda(n_j)}{\sqrt{n_j}} \widehat f(\frac{\log n_j}{\log q}) \\
+ \O\left(\frac 1{q}  \left(\frac 1{\log q} \sum_n
\frac{\Lambda(n)}{\sqrt{n}} \widehat f(\frac{\log n}{\log q})
\right)^m \right) \;.
\end{multline}

To bound the remainder term, use $\supp \widehat
f\subseteq[-\alpha, \alpha]$ with $\alpha<2/m$ to estimate the sum
$$
\sum_n  \frac{\Lambda(n)}{\sqrt{n}} \widehat f(\frac{\log n}{\log
q}) \ll \sum_{n\ll q^\alpha} \frac{\Lambda(n)}{\sqrt{n}} \ll
q^{\alpha/2}
$$
Thus the $\O$-term in \eqref{J(S)+O} is bounded by
$$
\frac 1{q(\log q)^m}  q^{m\alpha/2}\ll q^{-(1-m\alpha/2)} \to 0
$$

Thus we find
$$
\aveq{J(S)} = (-\frac 1{\log q})^m \sum_{\prod_{j\in S}n_j\equiv
\prod_{i\notin S}n_i\mod q} \prod_{j=1}^m
\frac{\Lambda(n_j)}{\sqrt{n_j}} \widehat f(\frac{\log n_j}{\log
q}) + \O(q^{-(1-m\alpha/2)})
$$

We split the sum for $\aveq{J(S)}$ into   terms $J_{eq}(S)$ where
we have equality $\prod_{j\in S}n_j= \prod_{i\notin S}n_i$ rather
then mere congruence modulo $q$, and the remaining terms:
$$
\aveq{J(S)} = J_{eq}(S) + J_{cong}(S) + \O(q^{-(1-m\alpha/2)})
$$
with
\begin{equation}\label{J_{eq}}
J_{eq}(S):= (-\frac 1{\log q})^m\sum_{ \prod_{j\in S}n_j=
\prod_{i\notin S}n_i} \prod_{j=1}^m
\frac{\Lambda(n_j)}{\sqrt{n_j}} \widehat f(\frac{\log n_j}{\log
q})
\end{equation}

\begin{equation*}
J_{cong}(S) := (-\frac 1{\log q})^m \sum_{\substack{ \prod_{j\in S}n_j\equiv \prod_{i\notin S}n_i\mod q\\
\prod_{j\in S}n_j\neq  \prod_{i\notin S}n_i}} \prod_{j=1}^m
\frac{\Lambda(n_j)}{\sqrt{n_j}} \widehat f(\frac{\log n_j}{\log
q})
\end{equation*}

\subsection{Eliminating congruential terms }
We will show that the  terms $J_{cong}(S)$ are negligible:
$$
J_{cong}(S)\ll q^{-(1-m\alpha/2)+\epsilon}
$$ for all $\epsilon>0$.

\begin{lem} \label{lem:non1}
Assume $XY=o(q^2)$. Then
$$
\sum_{\substack{ M<X,N<Y\\M\neq N\\M\equiv N\mod q }} \frac
1{\sqrt{MN}} \ll \frac{\sqrt{XY}}q
$$
and moreover if $\max(X,Y)=o(q)$ then the sum is empty, hence
equals zero.
\end{lem}
\begin{proof}
We may assume $X\leq Y$ and so certainly $X=o(q)$. Then in the sum
we must have $M<N$  since otherwise, $M=kq+N$ with $k\geq 1$ and
so $q<M<X=o(q)$ which gives a contradiction. If $Y=o(q)$ then
likewise the sum is empty. Thus we now assume that $Y\gg q$. In
that case we write $N=kq+M$, $1\leq k \leq Y/q$. Then our sum is
\begin{equation*}
\begin{split}
\sum_{M<X} \frac 1{\sqrt{M}} \sum_{k\leq Y/q} \frac 1{\sqrt{kq+M}}
&\ll \sum_{M<X} \frac 1{\sqrt{M}}  \sum_{k\leq Y/q} \frac
1{\sqrt{kq}}
\\
& \ll \sum_{M<X} \frac 1{\sqrt{M}} \frac{\sqrt{Y}}{q} \ll
\frac{\sqrt{XY}}q
\end{split}
\end{equation*}

as required.
\end{proof}

\begin{lem} Assume that $r+s=m$, and that $\alpha< 2/m$. Then
$$
\frac 1{(\log q)^m}
\sum_{\substack{ n_1,\dots,n_r,m_1,\dots m_s <q^\alpha \\
\prod n_i \equiv \prod m_j \mod q\\ \prod n_i \neq \prod m_j}}
\prod \frac{\Lambda(n_i)}{\sqrt{n_i}}
\prod\frac{\Lambda(m_j)}{\sqrt{m_j}} \ll q^{ m\alpha/2
-1+\epsilon}
$$
\end{lem}
\begin{proof}
We replace the sum over {\em prime powers} by the sum over all
{\em integers} and since all variables are bounded by $q^\alpha$,
we replace $\Lambda(n)$ by $\log q$. Thus our sum is $\ll$ than
$$
\sum_{\substack{ n_1,\dots,n_r,m_1,\dots m_s <q^\alpha \\
\prod n_i \equiv \prod m_j \mod q\\ \prod n_i \neq \prod m_j}}
\prod \frac{1}{\sqrt{n_i}} \prod\frac{1}{\sqrt{m_j}} \;.
$$
Now set $N=\prod n_i$, $M=\prod m_j$ and sum separately over those
tuples $n_1,\dots ,n_r$ with $N,M$ fixed. The number of such
tuples is $\ll q^\epsilon$ for all $\epsilon>0$. Thus our sum is
$\ll$ than
$$
q^\epsilon \sum_{\substack{M <q^{s\alpha}, N<q^{r\alpha} \\M\equiv
N \mod q\\M\neq N }} \frac 1{\sqrt{MN}}
$$
which by Lemma~\ref{lem:non1} is $\ll q^{\epsilon+ m\alpha/2 -1}$
for all $\epsilon>0$, since $q^{r\alpha}\cdot q^{s\alpha}=o(q^2)$
if $m\alpha = (r+s)\alpha<2$.
\end{proof}

Thus we find that
\begin{equation}\label{eq:Wf = J_eq + error}
\aveq{ (\Wf^{osc} )^m } = \sum_{S\subset\{1,\dots ,m\}} J_{eq}(S)
+ \O( q^{-(1-m\alpha/2)+\epsilon})
\end{equation}
for all $\epsilon>0$.

\subsection{Reduction to diagonal terms}

Fix a subset $S\subseteq\{1,\dots,m\}$. The sum in $J_{eq}(S)$
\eqref{J_{eq}} is over tuples $(n_1,\dots, n_m)$ which satisfy
$\prod_{j\in S} n_j = \prod_{i\notin S} n_i$. We say that there is
a {\em perfect matching} of terms if there  is a bijection
$\sigma$ of $S$ onto its complement $S^c$ in $\{1,\dots, m\}$ so
that $n_j=n_{\sigma(j)}$, for all $j\in S$. This can happen only
if $m=2k$ is {\em even} and $\# S = \# S^c=k$.

Decompose
\begin{equation}\label{Decomposing J_{eq}(S)}
J_{eq}(S)=J_{diag}(S)+J_{non}(S)
\end{equation}
where $J_{diag}(S)$ is the sum of matching terms - the diagonal
part of the sum (nonexistent for most $S$), and $J_{non}(S)$ is
the sum over the remaining, nonmatching, terms.

\subsection{Diagonal terms}
Assume that $m=2k$ is even. The diagonal terms  are the sum over
all $\binom{2k}{k}$ subsets $S\subset \{1,\dots 2k\}$ of
cardinality $k=m/2$ and for each such subset $S$, $J_{diag}(S)$ is
the sum over all $k!$ bijections $\sigma:S\to S^c$ of $S$ onto its
complement, of terms
$$
\left( \frac 1{(\log q)^2} \sum_n \frac {\Lambda(n)^2}{n} \widehat
f(\frac {\log n}{\log q})^2 \right)^k
$$
We evaluate each factor by using the Prime Number Theorem:
\begin{equation}\label{matching terms}
\begin{split}
\frac 1{(\log q)^2} \sum_n \frac {\Lambda(n)^2}{n} \widehat
f(\frac {\log n}{\log q})^2
&\sim\frac 1{(\log q)^2} \int_2^\infty \frac{\log t}{t} \widehat f(\frac {\log t}{\log q})^2 \;\d t\\
&\sim\int_0^\infty u\widehat f(u)^2 \;\d u
\end{split}
\end{equation}
Since our function is even and supported inside
$(-2/m,2/m)\subseteq(-1,1)$ (since $m\geq 2$), we can  rewrite
this as
$$
\frac 12 \int_{-\infty}^\infty \min(1,|u|) \widehat f(u)^2 \;\d
u=:\sigma(f)^2/2
$$
This shows that for $m=2k$ even we have as $q\to \infty$  that
$$
\sum_{S\subset\{1,\dots ,m\}} J_{diag}(S) \to  \frac{(2k)!}{2^k
k!} \sigma(f)^{2k }
$$

Below we will show that the nondiagonal terms $J_{non}(S)$ are
negligible, and hence by \eqref{eq:Wf = J_eq + error} and
\eqref{Decomposing J_{eq}(S)} we will have thus proved
Theorem~\ref{thm-L:Gaussian moments}.
\end{proof}

\subsection{Bounding the off-diagonal terms $J_{non}(S)$}

We will show that
\begin{lem}\label{bounding J_{non}}
$$
J_{non}(S) \ll \frac 1{\log q}
$$
\end{lem}
\begin{proof}

Since
$$
\frac 1{\log q} \sum_{p}\sum_{ k\geq 3} \frac{\log p}{p^{k/2}} \ll
\frac 1{\log q} \sum_p \frac {\log p}{ p^{3/2}} \ll \frac 1{\log
q}
$$
the contribution of cubes and higher prime powers to
\eqref{J_{eq}} is negligible, and we may assume in $J_{non}(S)$
that the $n_i$ are either prime or squares of primes (upto a
remainder of $\O(1/\log q)$). By the Fundamental Theorem of
Arithmetic, an equality $\prod_{i\in S}n_i = \prod_{j\in S^c}n_j$
forces some of the terms to match, and unless there is a perfect
matching of all terms, the remaining integers satisfy equalities
of the form $n_1n_2=n_3$ with $n_1=n_2=p$ prime and  $n_3=p^2$  a
square of that prime. Thus upto a remainder of $\O(1/\log q)$,
$J_{non}(S)$ is a sum of terms of the form
$$
\left( \frac 1{(\log q)^2} \sum_{\substack{p\\k=1,2}} \frac {(\log
p)^2}{p^k} \widehat f(\frac{\log p^k}{\log q})^2 \right)^u \cdot
\left( \frac 1{(\log q)^3} \sum_p \frac {(\log p)^3}{p^2} \widehat
f(\frac{\log p}{\log q})^2 \widehat f(\frac{\log p^2}{\log q})
\right)^v
$$
with $2u+3v=m$, and $v\geq 1$.

We showed \eqref{matching terms} that the matching terms have an
asymptotic value, hence are bounded. We bound the second type of
term  by
$$
\frac 1{(\log q)^3} \sum_p \frac {(\log p)^3}{p^2} \widehat
f(\frac{\log p}{\log q})^2 \widehat f(\frac{\log p^2}{\log q}) \ll
\frac 1{(\log q)^3} \sum_p \frac {(\log p)^3}{p^2} \ll \frac
1{(\log q)^3}
$$
Thus as long as $v\geq 1$  (that is if there is no perfect
matching of all terms), we get that the contribution of
$J_{non}(S)$ is $\O(1/\log q)$. This proves Lemma~\ref{bounding
J_{non}}.
\end{proof}

\section{The random matrix model}\label{sect:rmt_model}

Let $f(x)$ be an even real function subject to the decay condition
that there exists a fixed $\epsilon>0$ and $A>0$ such that
\begin{equation}\label{eq:f_decay}
f(x) < A (1+|x|)^{-(1+\epsilon)} \text{ for all } x \in \R
\end{equation}

Define
\begin{equation*}
F_N(\t) := \sum_{j=-\infty}^\infty f\left(\tfrac{N}{2\pi} (\t+2\pi
j)\right)
\end{equation*}
so that $F_N(\t)$ is $2\pi$-periodic. Define
\begin{equation*}
\Zf(U) := \sum_{j=1}^N F_N(\t_j)
\end{equation*}
where  $U$ is an $N\times N$ unitary matrix with eigenangles
$\t_1,\dots,\t_N$. This is the random matrix equivalent of
$\Wf(\chi)$.

The Fourier coefficients of $F_N(\t)$ are
\begin{align*}
a_{n,N} &= \frac{1}{2\pi} \int_{-\pi}^\pi F_N(\t) e^{-\i n\t} \;\d \t\\
&= \frac{1}{N} \int_{-\infty}^{\infty} f(x) e^{-2\pi\i n x /N} \;\d x\\
&= \frac{1}{N} \^ f\left(\frac{n}{N}\right)
\end{align*}
and so, if the matrix $U$ has eigenangles $\t_1,\dots,\t_N$,
\begin{align}
\Zf(U) &:= \sum_{j=1}^N F_N(\t_j)\notag\\
&= \sum_{n=-\infty}^\infty \frac{1}{N} \^
f\left(\frac{n}{N}\right) \Tr U^n .\label{eq:Zf_fourier}
\end{align}

Since
\begin{equation*}
\E\left\{\Tr U^n\right\} =
\begin{cases}
 N, & n=0\\
 0, & \text{otherwise}
\end{cases}
\end{equation*}
we have thus proven
\begin{thm}\label{thm:Zf_mean}
\begin{equation*}
\E \left\{ \Zf \right\} = \^ f(0)
\end{equation*}
\end{thm}
The definition of Fourier transform we use is such that $\^f(0) =
\int_{-\infty}^\infty f(x)\;\d x$, so this Theorem is in perfect
agreement with Theorem \ref{thm:expectation}.

\begin{thm}\label{thm:Zf_var}
The variance of $\Zf$ tends to $\sigma^2$ as $N\to\infty$, where
\begin{equation}
\sigma^2 = \int_{-\infty}^\infty \min(|u|,1) \^f(u)^2\;\d u
\label{eq:RMT_sigma^2}
\end{equation}
\end{thm}
\begin{proof}
Since \cite{DE,Haake,Rains}
\begin{equation*}
\E\left\{\Tr U^n \Tr U^m\right\} =
\begin{cases}
N^2 & \text{ if } n=m=0\\
|n| & \text{ if } n=-m \text{ and } |n|\leq N\\
N & \text{ if } n=-m \text{ and } |n|\geq N\\
0 & \text{ otherwise}
\end{cases}
\end{equation*}
we have
\begin{align*}
\E\left\{ \left(\Zf-\^f(0)\right)^2\right\}&= \sum_{\substack{n=-\infty\\n\neq 0}}^\infty \frac{1}{N^2} \^f\left(\frac{n}{N}\right) \^f\left(-\frac{n}{N}\right) \min(|n|,N)\\
&\to \int_{-\infty}^\infty \min(|u|,1) \^f(u)^2\;\d u
\end{align*}
the last line following from the definition of a Riemann integral,
and from the fact that $f(x)$ is even.
\end{proof}
Note that this is the same as the variance of $\Wf(\chi)$ (Theorem
\ref{thm:Wf_variance}) when $\^f$ is restricted to have support
contained in $[-1,1]$.

\subsection{Mock-Gaussian behaviour}\label{sect:Zf_mock_Gaussian}

 From \eqref{eq:Zf_fourier} and Theorem \ref{thm:Zf_mean}, the
$m$th centered moment is
\begin{multline*}
\E\left\{\left( \Zf - \E \{\Zf\} \right)^m \right\}\\
=\sum_{\substack{n_1=-\infty\\n_1\neq 0}}^\infty \dots
\sum_{\substack{n_m=-\infty\\n_m\neq 0}}^\infty
\frac{1}{N}\^f\left(\frac{n_1}{N}\right)\dots\frac{1}{N}\^f\left(\frac{n_m}{N}\right)
\E\left\{\Tr U^{n_1} \dots \Tr U^{n_m} \right\}
\end{multline*}
The following two Lemmas will enable us to calculate these
moments, under certain restrictions on the support of $\^f(u)$.

\begin{lem}\label{lem:sum_n_zero}
If $\sum_{j=1}^m n_j \neq 0$ then
\begin{equation*}
\E \left\{ \prod_{j=1}^m  \Tr U^{n_j} \right\} = 0
\end{equation*}
\end{lem}
\begin{proof}
By rotation invariance of Haar measure, the left hand-side is left
unchanged by multiplication by $I e^{\i\t}$ (where $\t$ is an
arbitrary angle, and $I$ is the identity matrix). Since $\Tr
\left((U I e^{\i\t})^n\right) = e^{\i n\t} \Tr U^n$, this means
$$
\E \left\{ \prod_{j=1}^m  \Tr U^{n_j} \right\} =
\exp\left(\i\t\sum_{j=1}^m n_j\right) \E \left\{ \prod_{j=1}^m
\Tr U^{n_j} \right\}
$$
which is true only if either both sides are zero, or if
$\sum_{j=1}^m n_j = 0$.
\end{proof}

\begin{lem}{\bf (Diaconis, Shahshahani \cite{DS,DE}).}\label{lem:DE}
For $a_j , b_j \in \{0,1,2,\dots\}$, if
\begin{equation*}
N \geq \max\left(\sum_{j=1}^k j a_j \ , \ \sum_{j=1}^k j
b_j\right) ,
\end{equation*}
then
\begin{equation*}
\E\left\{ \prod_{j=1}^k \left(\Tr U^{j} \right)^{a_j} \left(\Tr
U^{-j} \right)^{b_j} \right\} = \delta_{a,b} \prod_{j=1}^k j^{a_j}
a_j!
\end{equation*}
where $\delta_{a,b}=1$ if $a_j=b_j$ for $j=1,\dots,k$, and
$\delta_{a,b}=0$ otherwise.
\end{lem}

\begin{thm}\label{thm:rmt_mock_Gaussian}
For any integer $m\geq 2$, if $\supp \^f(u) \subseteq [-2/m,2/m]$,
then
\begin{equation*}
\lim_{N\to\infty} \E \left\{ \left( \Zf-\^ f(0) \right)^m \right\}
=
\begin{cases}
0 & \text{ if } m \text{ odd}\\
\frac{(2k)!}{2^{k}  k!} \sigma^m & \text{ if } m=2k, k\geq 1
\text{ an integer}
\end{cases}
\end{equation*}
where $\sigma^2$, the variance, is given by
\eqref{eq:RMT_sigma^2}.
\end{thm}

\begin{proof}
The restriction on the support of $\^ f(u)$ gives
\begin{multline}\label{eq:m_moment_expand}
\E \left\{ \left( \Zf-\^ f(0) \right)^m \right\} \\
= \frac{1}{N^m} \sum_{\substack{n_1=-2N/m\\n_1\neq 0}}^{2N/m}
\dots \sum_{\substack{n_m=-2N/m\\n_m\neq 0}}^{2N/m}
\^f\left(\frac{n_1}{N}\right) \dots \^f\left(\frac{n_m}{N}\right)
\E \left\{ \Tr U^{n_1} \dots \Tr U^{n_m}\right\}
\end{multline}

Lemma~\ref{lem:sum_n_zero} means that to have a non-zero
contribution, $\sum n_j=0$, and so
\begin{equation*}
\max_{\substack{n_j \\|n_j|<2N/m\\\sum n_j = 0}}
\left\{\sum_{j=1}^m n_j \I_{\{n_j>0\}}\right\} \leq \frac{m}{2}
\frac{2N}{m} = N
\end{equation*}
the maximum is obtained by all the positive terms equal to $2N/m$,
and all the negative terms to $-2N/m$. (This maximum is obtainable
only if $m$ is even). Thus we see that the support restriction
means all the nonzero terms in \eqref{eq:m_moment_expand} can be
calculated using Lemma~\ref{lem:DE}.

To obtain anything nonzero using Lemma~\ref{lem:DE}, there must be
a bijection $\sigma$ mapping $\{1,\dots,m\}$ into itself so that
$n_j=-n_{\sigma(j)}$ for all $j$. Note that no $n_j$ can equal
zero, since this is expressly forbidden in
\eqref{eq:m_moment_expand}.

For odd $m$, it is impossible to pair off the $n_j$ without having
at least one $n_j=0$. Therefore \eqref{eq:m_moment_expand} is zero
for $m$ odd.

For even $m=2k$, assume the $n_j$ are such that they can be paired
off, and relabel so that $r_1=n_{j_1}$ where $j_1$ is the smallest
number such that $n_{j_1}>0$, $r_2=n_{j_2}$ where $j_2$ is the
second smallest number such that $n_{j_2}>0$ etc. There are
$\binom{2k}{k}$ ways of arranging the positive $n_j>0$ to give the
same $r_i$. The number of ways of ordering the negative $n_j$ such
that each positive term has a negative partner equals
$\frac{k!}{b_1! \ b_2! \dots}$ where $b_i = \#\{j \ : \ n_j =
-i\}$. Therefore, after reordering, \eqref{eq:m_moment_expand}
equals
\begin{align*}
\sum_{r_1=1}^{N/k} \dots \sum_{r_k=1}^{N/k} \binom{2k}{k} \frac{k!}{b_1! \ b_2! \dots} &\E\left\{\left|\Tr U^{r_1}\right|^2 \dots\left| \Tr U^{r_k}\right|^2\right\} \prod_{i=1}^k \frac{1}{N^2} \^f(\frac{r_i}{N}) \^f(\frac{-r_i}{N})\\
&= \frac{(2k)!}{k!} \left( \sum_{r=1}^{N/k}  r \left|\frac{1}{N}\^f\left(\frac{r}{N}\right)\right|^2\right)^k\\
&\to \frac{(2k)!}{k! 2^k} \sigma^2
\end{align*}
since
\begin{equation*}
\frac{1}{b_1! \ b_2! \dots} \E\left\{\left|\Tr U^{r_1}\right|^2
\dots\left| \Tr U^{r_k}\right|^2\right\} = \prod_{j=1}^k r_j
\end{equation*}
by Lemma~\ref{lem:DE}, and since
\begin{equation*}
\sum_{r=1}^{N/k} r
\left|\frac{1}{N}\^f\left(\frac{r}{N}\right)\right|^2 \sim
\frac{1}{2}\int_{-1/k}^{1/k} |u| \left|\^ f(u)\right|^2 \;\d u =
\tfrac{1}{2}\sigma^2
\end{equation*}
when $\supp \^f \subseteq [1/k,1/k]$, where $\sigma^2$ is given by
\eqref{eq:RMT_sigma^2}.
\end{proof}

\begin{remark}
One can also prove Theorem~\ref{thm:rmt_mock_Gaussian} by a
completely different method, using techniques found in
\cite{Sosh}; we need to take this route when dealing with the
other classical groups in \cite{hr2}.
\end{remark}

\section{Unrestricted moments of $\Zf(U)$}\label{sect:Z_f_not_fourier}

In this section we will calculate the (uncentered) $m$th moment of
$\Zf(U)$ without restriction on the support. This allows us to
conjecture an extension to Theorems \ref{thm:expectation} and
\ref{thm-L:Gaussian moments}. In particular, it appears that the
$m$th centered moment of $\Wf(\chi)$ is not Gaussian outside of
the range given in Theorem \ref{thm-L:Gaussian moments}, which
would imply that $\Wf(\chi)$ does not converge to a normal
distribution.

We wish to calculate the (uncentered) $m$th moment of $\Zf(U)$.
\begin{align}
M_m &:= \lim_{N\to\infty} \E \left\{ \left( \Zf \right)^m \right\}\notag\\
&= \lim_{N\to\infty} \E \left\{ \sum_{i_1=1}^N \dots
\sum_{i_n=1}^N
 F_N(\t_{i_1}) \dots F_N(\t_{i_n}) \right\}  \label{eq:M_m_not_Fourier}
\end{align}

To evaluate $M_m$, we use the $r$-point correlation function of
Dyson:
\begin{lem}\label{lem:Dyson}
{\bf (Dyson).} For an arbitrary function $g$ of $r$ variables
which is $2\pi$--periodic in all its variables,
\begin{multline*}
\E\left\{\sum_{\substack{i_1,\dots,i_r=1\\i_j \text{ all distinct}}}^N g(\t_{i_1},\dots,\t_{i_r}) \right\} \\
=\frac{1}{(2\pi)^r} \int_{-\pi}^\pi \dots\int_{-\pi}^\pi
g(\t_1,\dots,\t_r) R_r^{(N)}(\t_1,\dots,\t_r) \;\d \t_1\dots\d
\t_r
\end{multline*}
where
\begin{equation*}
R_r^{(N)}(\t_1,\dots,\t_r) = \det\left\{ S_N(\t_j-\t_i)
\right\}_{1\leq i,j\leq r}
\end{equation*}
with
\begin{equation*}
S_N(x) = \frac{\sin(Nx/2)}{\sin(x/2)}
\end{equation*}
\end{lem}

Note that the sums in \eqref{eq:M_m_not_Fourier} range
unrestrictedly over all variables (they include both diagonal and
off-diagonal terms), whereas Lemma \ref{lem:Dyson} requires the
sums to be over distinct variables (off-diagonals only). We
overcome this problem by summing over the diagonals separately.

\begin{defn}\label{defn:partition}
$\sigma$ is said to be a set partition of $m$ elements into $r$
non-empty blocks if
\begin{equation*}
\sigma : \{1,\dots,m\} \longrightarrow \{1,\dots,r\}
\end{equation*}
satisfying
\begin{enumerate}
\item For every $q\in\{1,\dots,r\}$ there exists at least one $j$ such
 that $\sigma(j)=q$ (this is the non-emptiness of the blocks).
\item For all $j$, either $\sigma(j)=1$ or there exists a $k<j$ such
 that $\sigma(j)=\sigma(k)+1$. (Roughly speaking, if we think of
 $\{1,\dots,r\}$ as denoting ordered pigeonholes, then $\sigma(j)$
 either goes into a non-empty pigeonhole, or into the next empty hole).
\end{enumerate}
The collection of all set partitions of $m$ elements into $r$
blocks is denoted $P(m,r)$.
\end{defn}

\begin{remark}
The number of $\sigma\in P(m,r)$  is equal to $S(m,r)$, a Stirling
number of the second kind. The number of set partitions of $m$
elements into any number of non-empty blocks is $\sum_{r=1}^m
S(m,r) = B_m$, a Bell number.
\end{remark}

\begin{lem}\label{lem:sum_distinct}
For any function $g$ of $m$ variables,
\begin{equation*}
\sum_{j_1,\dots,j_m} g(x_{j_1},\dots,x_{j_m}) = \sum_{r=1}^m
\sum_{\sigma \in P(m,r)} \sum_{\substack{\i_1,\dots,i_r\\i_j
\text{ all distinct}}}
g(x_{i_{\sigma(1)}},\dots,x_{i_{\sigma(m)}})
\end{equation*}
\end{lem}
\begin{proof}
Each term on the LHS appears once and only once on the RHS, so
they are equal.
\end{proof}

\begin{thm}\label{thm:M_m}
\begin{equation*}
M_m = \sum_{r=1}^m \nidotsint{r}_{-\infty}^\infty
R_r(x_1,\dots,x_r) \sum_{\sigma \in P(m,r)} \prod_{q=1}^r
f^{\lambda_q}(x_{q}) \;\d x_{q}
\end{equation*}
where $\l_q = \#\{ j : \sigma(j)=q \}$, and where
\begin{equation*}
R_r(x_1,\dots,x_r) = \det\left\{
\frac{\sin(\pi(x_j-x_i))}{\pi(x_j-x_i)} \right\}_{1\leq i,j\leq r}
\end{equation*}
\end{thm}

\begin{proof}
Recall \eqref{eq:M_m_not_Fourier}, that
\begin{equation*}
M_m = \lim_{N\to\infty} \E \left\{ \sum_{i_1=1}^N \dots
\sum_{i_m=1}^N F_N(\t_{i_1}) \dots F_N(\t_{i_m}) \right\}
\end{equation*}

Lemma \ref{lem:sum_distinct} gives
\begin{multline} \label{eq:expectation_F_N_distinct}
\E \left\{ \sum_{i_1=1}^N \dots \sum_{i_m=1}^N F_N(\t_{i_1}) \dots F_N(\t_{i_m}) \right\} \\
= \sum_{r=1}^m \sum_{\sigma\in P(m,r)} \E\left\{ \sum_{\substack{i_1,\dots,i_r=1\\i_j \text{ all distinct}}}^N F_N(\t_{i_{\sigma(1)}}) \dots F_N(\t_{i_{\sigma(m)}}) \right\}\\
= \sum_{r=1}^m \sum_{\sigma\in P(m,r)} \E\left\{
\sum_{\substack{i_1,\dots,i_r=1\\i_j \text{ all distinct}}}^N
F_N^{\lambda_1}(\t_{i_1}) \dots F_N^{\lambda_r}(\t_{i_r}) \right\}
\end{multline}
where $\lambda_q = \#\{ j : \sigma(j)=q \}$.  Lemma
\ref{lem:Dyson} now applies, and gives
\begin{multline*}
\E\left\{ \sum_{\substack{i_1,\dots,i_r=1\\i_j \text{ all distinct}}}^N F_N^{\lambda_1}(\t_{i_1}) \dots F_N^{\lambda_r}(\t_{i_r}) \right\}\\
= \frac{1}{(2\pi)^r}\idotsint_{-\pi}^\pi F_N^{\lambda_1}(\t_{1}) \dots F_N^{\lambda_r}(\t_{r})  R_r^{(N)}(\t_1,\dots,\t_r) \;\d\t_1\dots\d\t_r\\
= \frac{1}{N^r} \idotsint_{-N/2}^{N/2} R_r^{(N)}\left(\frac{2\pi
x_1}{N},\dots,\frac{2\pi x_r}{N}\right) \prod_{q=1}^r
F_N^{\lambda_q}\left(\frac{2\pi x_{q}}{N}\right) \;\d x_q
\end{multline*}
upon change variables to $x_n = \frac{N}{2\pi}\t_n$. Now,
\begin{align*}
F_N\left(\frac{2\pi x}{N}\right) &= \sum_{j=-\infty}^\infty f(x+Nj)\\
&= f(x) + \O\left(\frac{1}{N^{1+\epsilon}}\right)
\end{align*}
uniformly for all $x\in[-N/2,N/2]$, due to the decay condition on
$f$, \eqref{eq:f_decay}.

Since
\begin{equation*}
\lim_{N\to\infty} \frac{1}{N^r} R_r^{(N)}\left(\frac{2\pi
x_1}{N},\dots,\frac{2\pi x_r}{N}\right) =
R_r\left(x_1,\dots,x_r\right)
\end{equation*}
where
\begin{equation*}
R_r(x_1,\dots,x_r) = \det\left\{
\frac{\sin(\pi(x_j-x_i))}{\pi(x_j-x_i)} \right\}_{1\leq i,j\leq r}
\end{equation*}
we have
\begin{multline}\label{eq:N_infty_expectation_F_N_distinct}
\lim_{N\to\infty} \E\left\{ \sum_{\substack{i_1,\dots,i_r=1\\i_j \text{ all distinct}}}^N F_N^{\lambda_1}(\t_{i_1}) \dots F_N^{\lambda_r}(\t_{i_r}) \right\} \\
= \idotsint_{-\infty}^{\infty} R_r\left(x_1,\dots,x_r\right)
\prod_{q=1}^r f^{\lambda_q}(x_q) \;\d x_q
\end{multline}

Hence, combining \eqref{eq:M_m_not_Fourier},
\eqref{eq:expectation_F_N_distinct} and
\eqref{eq:N_infty_expectation_F_N_distinct}
\begin{align}\notag
M_m &= \lim_{N\to\infty} \E \left\{ \sum_{i_1=1}^N \dots \sum_{i_m=1}^N F_N(\t_{i_1}) \dots F_N(\t_{i_m}) \right\}\\
&= \sum_{r=1}^m  \sum_{\sigma \in P(m,r)}
\nidotsint{r}_{-\infty}^\infty  R_r(x_1,\dots,x_r) \prod_{q=1}^r
f^{\lambda_q}(x_{q}) \;\d x_{q} \label{eq:M_m_answer_not_Fourier}
\end{align}
as required.
\end{proof}

\begin{remark}
One can show that the moments $M_m$ of $\Zf$ uniquely determine
the distribution $\Zf$ weakly converges to as $N\to\infty$.
\end{remark}

\section{Application: small first zeros of $L(s,\chi)$}\label{sect:small_zeros}

In this section we will apply the results of \S\S\ref{sect:W_mean}
and \ref{sect:variance of Wf} to show that, under the assumption
of GRH there exist Dirichlet $L$--functions whose first zero is
lower than the expected height. Small gaps between high zeros of
the Riemann zeta function (which also obey unitary statistics)
have been much studied. Montgomery \cite{Mont73} showed that an
infinite number of zeros are less than $0.68$ times their average
spacing. This was reduced to $0.5179$ by Montgomery and Odlyzko
\cite{MontOd}; to $0.5171$ by Conrey, Ghosh and Gonek \cite{CGG};
and to $0.5169$ by Conrey and Iwaniec, as announced in \cite{CI}.
Conrey, Ghosh, Goldston, Gonek and Heath-Brown \cite{CGGGH-B}
showed that a positive proportion of zeros are less than $0.77$
times the average spacing, a result improved to $0.6878$ by
Soundararajan \cite{Sound}. We should perhaps point out that the
main difference between gaps between the zeta zeros, and the
height of the lowest Dirichlet zero is that the point $1/2$ in not
expected to ``repel'' low-lying zeros.

\subsection{Infinitely many small first zeros}

Using Theorem \ref{thm:expectation} we are able to obtain some
partial results for extreme low-lying zeros of Dirichlet
$L$--functions.

\begin{thm}\label{thm:inf_less_1/4R}
Assume GRH. If
\begin{equation}\label{eq:average of Wf assumption}
\lim_{q\to\infty} \ave{\Wf}_q = \int_{-\infty}^\infty f(x) \;\d x
\end{equation}
for all admissible functions $f$ with $\supp \^f \subseteq
[-2R,2R]$, then
\begin{equation*}
\liminf_{q\to\infty} \min_{\chi\neq\chi_0} x_{\chi,1} \leq
\frac{1}{4R}
\end{equation*}
where the minimum of the first zero of $L(s,\chi)$ is taken over
all non-trivial characters modulo $q$.
\end{thm}

\begin{proof}
Let $\^g(u)$ be an even, continuous function, with $\supp \^g
\subset [-R,R]$, and such that $\^g(u)$ is differentiable in
$[-R,R]$ and $g(x) \ll |x|^{-3/2-\delta}$, $\delta>0$.

Let
\begin{equation*}
B := \sqrt{\frac{ \int_0^\infty x^2 g^2(x) \;\d x}{\int_0^\infty
g^2(x) \;\d x }} = \sqrt{\frac{\frac{1}{4\pi^2} \int_0^\infty
\^g'(u)^2\;\d u}{\int_0^\infty \^g(u)^2 \;\d u}}
\end{equation*}
so, by assumptions on $\^g(u)$ and its derivative, we see that
$\b$ is a strictly positive finite real number.

Define, for $\beta>B$,
\begin{equation*}
f(x) = (x^2-\b^2)g^2(x)
\end{equation*}
so that $f$ has the properties
\begin{equation}\label{eq:int f < 0}
\int_0^\infty f(x) \;\d x = -(\beta^2-B^2)\int_0^\infty g^2(x)\;\d
x < 0
\end{equation}
and
\begin{equation*}
f(x)
\begin{array}{c}
 \leq\\
 \geq
\end{array}
0 \text{ for } |x|
\begin{array}{c}
\leq\\
\geq
\end{array}
\b
\end{equation*}
Note that the conditions on $\^g(u)$ mean that $f$ is an
admissible function.

Observe that
\begin{align*}
\^f(u) &= \frac{-1}{4\pi^2} \frac{\d^2}{\d u^2} (\^g\star \^g)(u) - \b^2 (\^g\star\^g)(u)\\
&= \tfrac{-1}{4\pi^2} (\^g'\star \^g')(u) - \b^2 (\^g\star\^g)(u)
\end{align*}
where $(\^g\star\^g)(u)$ is the convolution of $\^g$ with itself.
Since differentiation and multiplication by a constant does not
increase the support of a function, we may conclude that $\supp
\^f \subset [2R,2R]$ since, by assumption, $\supp \^g \subset
[-R,R]$.

Therefore, by assumption \eqref{eq:average of Wf assumption} and
by \eqref{eq:int f < 0},
\begin{equation*}
\frac{1}{q-2} \sum_{\chi\neq\chi_0} \sum_{j\geq 1}
f\left(x_{\chi,j}\right) \sim \int_0^\infty f(x) \;\d x <0 .
\end{equation*}
By the assumption of GRH all the $x_{\chi,j}$ are real, and so we
may conclude that there exists a $q_0$ such that for all $q>q_0$,
\begin{equation*}
\frac{1}{q-2} \sum_{\chi\neq\chi_0} \sum_{\substack{j\geq
1\\x_{\chi,j}\leq\b}} f\left(x_{\chi,j}\right) < \frac{1}{q-2}
\sum_{\chi\neq\chi_0} \sum_{j\geq 1} f\left(x_{\chi,j}\right) < 0
\end{equation*}
and so, for all $q>q_0$ at least one $\chi$, a non-trivial
character modulo $q$, exists with $x_{\chi,1} \leq \beta$. (Note
that this method produces a non-vacuous result only if $\beta<1$,
since by definition, $\ave{x_{\chi,1} }_q \to 1$.) The Theorem
will follow if we can construct a $g(x)$ satisfying all the
conditions such that $B=1/4R$, by letting $\beta\to B$.

Taking $\^g(u) = \cos\left(\frac{\pi u}{2R}\right) \I_{\{ |u|\leq
R\} }$, so that
\begin{equation*}
g(x) = \frac{-4R\cos(2\pi xR)}{\pi(16x^2R^2-1)} ,
\end{equation*}
we see that
\begin{equation*}
B^2 = \frac{ \int_0^\infty x^2 g^2(x) \;\d x}{ \int_0^\infty
g^2(x) \;\d x} = \frac{1}{16 R^2}
\end{equation*}
This concludes the proof of Theorem \ref{thm:inf_less_1/4R}.
\end{proof}

\begin{remark}
Our choice of $\^g(u) = \cos\left(\frac{\pi u}{2R}\right) \I_{\{
|u|\leq R\} }$ was not an arbitrary one, as this is the optimizing
function for this method.
\end{remark}

\begin{cor}\label{thm:inf_less_1/4}
If the Generalised Riemann Hypothesis holds, then
\begin{equation*}
\liminf_{q\to\infty} \min_{\chi\neq\chi_0} x_{\chi,1} \leq
\frac{1}{4}
\end{equation*}
where the minimum of the first scaled zero of $L(s,\chi)$ is taken
over all non-trivial characters modulo $q$.
\end{cor}

\begin{proof}
By Theorem \ref{thm:expectation} we may take $R=1$ in Theorem
\ref{thm:inf_less_1/4R}.
\end{proof}

\begin{remark}
Random matrix theory suggests that
\begin{equation*}
\liminf_{q\to\infty} \min_{\chi\neq\chi_0} x_{\chi,1} = 0
\end{equation*}
\end{remark}

\subsection{Positive proportion of small first zeros}

Theorem \ref{thm:inf_less_1/4R} combined with Theorem
\ref{thm-L:Gaussian moments} allows us to deduce a statement about
a positive proportion (rather than just infinitely many) of the
$\chi$ have smaller than expected first zeros.

\begin{thm}\label{thm:pos_prop_small_zeros}
Assume GRH. For $\b\geq 0.633$,
\begin{multline*}
\liminf_{q\to\infty} \frac{1}{q-2} \# \left\{ \chi\neq\chi_0 \ : x_{\chi,1} < \beta \right\} \\
\geq
\frac{11\pi^2-3-72\b^2-88\pi^2\b^2-48\b^4+176\pi^2\b^4}{12\pi^2(4\b^2-1)^2}
\end{multline*}
\end{thm}
\begin{remark}
Random matrix theory suggests that a positive proportion of the
$\chi$ have $x_{\chi,1}<\beta$ for any $\beta>0$.
\end{remark}

\begin{proof}[Proof of Theorem \ref{thm:pos_prop_small_zeros}]
Take
\begin{equation*}
f_\b(x) = (x^2-\b^2)g^2(x)
\end{equation*}
where
\begin{equation*}
\^g(u) = \cos(\pi u) \I_{ \{ |u|\leq 1/2 \} }
\end{equation*}
(so $\^ f(u)$ has support in $[1,1]$, and $f(x)\leq 0$ for
$|x|\leq\b$, and $f(x)\geq 0$ otherwise).

As in the proof of Theorem \ref{thm:inf_less_1/4R} we have
\begin{align*}
\lim_{q\to\infty} \ave{\Wf}_q &= \int_{-\infty}^\infty f(x)\;\d x \\
&< 0 \text{ for $\b>1/2$.}
\end{align*}

By Theorem \ref{thm-L:Gaussian moments},
\begin{equation*}
\lim_{q\to\infty} \ave{ \left( \Wf - \ave{\Wf}_q \right)^{2} }_q =
\int_{-1}^1 |u| \left| \^f(u)\right|^2 \;\d u
\end{equation*}

Chebyshev's inequality gives
\begin{equation*}
\limsup_{q\to\infty} \frac{1}{q-2} \# \left\{ \chi\neq\chi_0 \ : \
\left|\Wf -  \ave{\Wf}_q\right| \geq \epsilon \right\} \leq \frac{
\int_{-1}^1 |u| \left| \^f(u)\right|^2 \;\d u}{\epsilon^2}
\end{equation*}
and so, using the fact that $f$ is even,
\begin{multline*}
\liminf_{q\to\infty} \frac{1}{q-2} \# \left\{ \chi\neq\chi_0 \ : \ \left|\sum_{j\geq 1}f(x_{\chi,j}) -  \int_0^\infty f(x)\;\d x \right| \leq \epsilon_1 \right\}\\
\geq 1- \frac{ \int_{-1}^1 |u| \left| \^f(u)\right|^2 \;\d
u}{4\epsilon_1^2}
\end{multline*}
where $\epsilon_1=\epsilon/2$.

If $\beta>1/2$, putting $\epsilon_1 = \left| \int_0^\infty
f(x)\;\d x\right|=\frac{1}{2}\left|\^f(0)\right|$, we get
\begin{multline*}
\liminf_{q\to\infty} \frac{1}{q-2} \# \left\{ \chi\neq\chi_0 \ : \ -2\left| \int_0^\infty f(x)\;\d x \right| \leq \sum_{j\geq 1}f(x_{\chi,j})  \leq 0 \right\}\\
\geq 1- \frac{ \int_{-1}^1 |u| \left| \^f(u)\right|^2 \;\d
u}{\^f(0)^2}
\end{multline*}
(Note we need GRH here, so that $\sum_{j\geq 1} f(x_{\chi,j})$ is
real). Since $f(x_{\chi,j})<0$ implies $x_{\chi,j}<\beta$ we may
conclude that, after working out the integrals on the right hand
side,
\begin{multline*}
\liminf_{q\to\infty} \frac{1}{q-2} \# \left\{ \chi\neq\chi_0 \ : x_{\chi,1} < \beta \right\} \\
\geq
1-\frac{3+\pi^2+72\b^2-8\pi^2\b^2+48\b^4+16\pi^2\b^4}{12\pi^2(4\b^2-1)^2}
\end{multline*}
The right hand side is greater than zero for
\begin{equation*}
\b \geq
\frac{1}{2}\frac{\sqrt{9+11\pi^2+2\sqrt{18+66\pi^2}}}{\sqrt{11\pi^2-3}}
\approx 0.633
\end{equation*}
as required.
\end{proof}

\begin{remark}
The test function we used in the proof is the optimum test
function for Theorem \ref{thm:inf_less_1/4R}, but that does not
necessarily make it the optimum test function here. Indeed, the
word ``optimum'' is not well defined here, as one can either try
to find a function that maximises the estimate of the proportion
of $\chi$ satisfying $x_{1,\chi}<\beta$, or one could try to find
a function which minimises the $\beta$ for which this method
proves a positive proportion of $x_{1,\chi}\leq \b$.
\end{remark}


\begin{thebibliography}{99}

\bibitem{CGG} J.B. Conrey, A. Ghosh and S.M. Gonek, {\em A note on
gaps between zeros of the zeta function}, Bull. London
Math. Soc. \textbf{16} (1984) 421--424 

\bibitem{CGGGH-B} J.B. Conrey, A. Ghosh, D. Goldston, S.M. Gonek and
D.R. Heath-Brown, {\em On the distribution of gaps between zeros of
the zeta-function}, Quart. J. Math. Oxford (2) \textbf{36} (1985)
43--51 

\bibitem{CI} J.B. Conrey and H. Iwaniec, {\em Spacing of zeros of
Hecke $L$--functions and the class number problem}, Acta
Arith. \textbf{103} (2002) no. 3, 259--312. 

\bibitem{DE} P. Diaconis and S.N. Evans, {\em Linear functionals of
eigenvalues of random matrices}, Trans. Amer. Math. Soc. 
\textbf{353} (2001) 2615--2633.  

\bibitem{DS} P. Diaconis and M. Shahshahani, {\em On the eigenvalues
of random matrices}, J. Appl. Probab. \textbf{31A} (1994)
49--62 

\bibitem{Haake} F. Haake, M. K\'us, H.-J. Sommers, H. Schomerus and
K. {\.Z}yczkowski  {\em Secular determinants of random unitary
matrices}, J. Phys. A \textbf{29} (1996) 
3641--3658.  

\bibitem{hr2} C.P. Hughes and Z. Rudnick, {\em Mock-Gaussian behaviour
for linear statistics of classical compact groups}, 
to appear in J. Phys. A. 


\bibitem{hr3} C.P. Hughes and Z. Rudnick, {\em Linear statistics for
zeros of Riemann's zeta function},  
C. R. Acad. Sci. Paris, Ser. {\bf I 335} (2002), 667-670.

\bibitem{ILS}
H. Iwaniec, W. Luo\ and\ P. Sarnak, {\em Low lying zeros of
families of $L$-functions}, Inst. Hautes \'Etudes Sci. Publ. Math.
No. \textbf{91} (2000), 55--131.


\bibitem{KatzSarnak} N.M. Katz and P. Sarnak, \textit{Random Matrices,
Frobenius Eigenvalues, and Monodromy}, (AMS Colloquium Publications,
1999). 

\bibitem{Mont73} H.L. Montgomery, {\em The pair correlation of zeros
of the zeta function}, Proc. Sym. Pure Math \textbf{24}
(1973) 181--193.  

\bibitem{MontOd} H.L. Montgomery and A.M. Odlyzko, ``Gaps between zeros of the zeta function'',  in \textit{Topics in Classical Number Theory}, Colloq. Math. Soc. Janos Bolyai \textbf{34} (1984) 1079--1106

\bibitem{MontVau} H.L. Montgomery and R.C. Vaughan, {\em The large
sieve}, Mathematika \textbf{20} (1973) 119--134. 

\bibitem{Oz} A. \"Ozl\"uck, {\em On the $q$-analogue of the pair
correlation conjecture}, J. Number Theory \textbf{59} (1996)
319--351.  

\bibitem{OzSn} A. \"Ozl\"uck and C. Snyder, {\em On the distribution
of the nontrivial zeros of quadratic $L$--functions close to the real
axis}, Acta Arith. \textbf{91} (1999) 209--228.  

\bibitem{Rains} E. Rains, {\em High powers of random elements of
compact Lie groups},  Probab. Theor. Rel. Fields \textbf{107} (1997)
219--241.  

\bibitem{Rubin} M. Rubinstein, {\em Low-lying zeros of $L$--functions
and random matrix theory}, Duke Math. J. \textbf{109} (2001) 147--181.   

\bibitem{selberg_L} A. Selberg, {\em Contributions to the theory of
Dirichlet's $L$--functions}, Skr. Norske Vid. Akad. Oslo. I. No. 3,
(1946) 1--62.    

\bibitem{Siegel} C.L. Siegel {\em The zeros of Dirichlet $L$-functions},
Annals of Math.  \textbf{46} (1945), 409--422.

\bibitem{Sound} K. Soundararajan, {\em On the distribution of gaps
between zeros of the Riemann zeta-function},
Quart. J. Math. Oxford (2) \textbf{47} (1996) 383--387.  

\bibitem{Sosh} A. Soshnikov, {\em Central limit theorem for local
linear statistics in classical compact groups and related
combinatorial identities}, Ann. Probab. \textbf{28} (2000) 1353--1370.     


\bibitem{Titchmarsh} E.C. Titchmarch {\em The zeros of Dirichlet's
$L$-functions}, Proc. London Math. Soc. \textbf{(2)} \textbf{32}
(1931), 488--500. 

\end{thebibliography}
\end{document}